# $U$-tests for variance components in one-way random effects models[*]

**Juvêncio S. Nobre**[1], **Julio M. Singer**[2] **and Mervyn J. Silvapulle**[3]

*Universidade Federal do Ceará, Universidade de São Paulo and Monash University*

**Abstract:** We consider a test for the hypothesis that the within-treatment variance component in a one-way random effects model is null. This test is based on a decomposition of a $U$-statistic. Its asymptotic null distribution is derived under the mild regularity condition that the second moment of the random effects and the fourth moment of the within-treatment errors are finite. Under the additional assumption that the fourth moment of the random effect is finite, we also derive the distribution of the proposed $U$-test statistic under a sequence of local alternative hypotheses. We report the results of a simulation study conducted to compare the performance of the $U$-test with that of the usual $F$-test. The main conclusions of the simulation study are that (i) under normality or under moderate degrees of imbalance in the design, the $F$-test behaves well when compared to the $U$-test, and (ii) when the distribution of the random effects and within-treatment errors are nonnormal, the $U$-test is preferable even when the number of treatments is small.

## 1. Introduction

Consider the one-way random effects model

(1.1) $$Y_{ij} = \mu + b_i + e_{ij}, \ i=1,\ldots k, \ j=1,\ldots,n_i \ (\geq 2)$$

with $b_i$ and $e_{ij}$, $i=1,\ldots,k$, $j=1,\ldots,n_i$, denoting independent random variables with null means and variances $\sigma_b^2$ and $\sigma_e^2$, respectively. Thus, $\mu$ is the mean response, $b_i$ represents the random effect of the $i$-th treatment, $e_{ij}$ represents a random measurement error associated with the $j$-th observation obtained under the $i$-th treatment, and $\sigma_b^2$ and $\sigma_e^2$ are between- and within-treatment variance components, respectively.

In general, data analysis based on such a model focuses on the estimation of $\mu$ and on testing the one-sided hypothesis (of no treatment effects)

(1.2) $$\mathcal{H}_0 : \sigma_b^2 = 0 \quad \text{vs} \quad \mathcal{H}_1 : \sigma_b^2 > 0.$$

Inference about variance components in random effects models, and more generally in linear mixed models, has a long history in the statistical literature. In this context, Searle et al. [21] provides an excellent overview of estimation and prediction while Khuri et al. [13] and Demidenko [6] present an extensive review of hypothesis testing; Silvapulle and Sen [36] consider a comprehensive treatment of testing under inequality constraints.

[*]Supported by Conselho Nacional de Desenvolvimento Científico e Tecnológico (CNPq) and Fundação de Amparo à Pesquisa do Estado de São Paulo (FAPESP), Brazil.
[1]Departamento de Estatística e Matemática Aplicada, Universidade Federal do Ceará, Brazil, e-mail: juvencio@ufc.br
[2]Departamento de Estatística, Universidade de São Paulo, Brazil, e-mail: jmsinger@ime.usp.br
[3]Department of Econometrics and Business Statistics, Monash University, P.O. Box 197, Caulfield East, Australia 3145, e-mail: Mervyn.Silvapulle@buseco.monash.edu.au

Under the additional assumption that $b_i$ and $e_{ij}$ follow normal distributions, the usual $F$-statistic for testing (1.2) is

$$(1.3) \qquad F = \frac{\mathrm{SQ}(b)/(k-1)}{\mathrm{SQ}(e)/(n-k)},$$

where $\mathrm{SQ}(b) = \sum_{i=1}^{k} n_i(\overline{Y}_{i.} - \overline{Y}_{..})^2$ and $\mathrm{SQ}(e) = \sum_{i=1}^{k}\sum_{j=1}^{n_i} Y_{ij}^2 - \sum_{i=1}^{k} n_i \overline{Y}_{i.}^2$ are, respectively, the between- and within-treatment sums of squares. The $F$-statistic (1.3) follows a central $F$ distribution with $k-1$ and $n-k$ degrees of freedom when $\mathcal{H}_0 : \sigma_b^2 = 0$ is true. In the balanced case (i.e., when $n_1 = \cdots = n_k$) the test is uniformly most powerful invariant (UMPI). This optimality property does not hold in the unbalanced case. Details may be found in Khuri et al. [13], for example.

We propose an alternative test based on the decomposition of $U$-statistics considered by Pinheiro et al. [19] in a nonparametric setup. Although it is not an exact test, it has good properties for moderate sample sizes and does not require the normality assumption. It is also valid for heteroskedastic random effects.

The class of $U$-statistics has its genesis in the papers of Halmos [11] and Hoeffding [12] and is well known for its simple structure and for the weak assumptions required for its use in statistical inference. It also provides a unified paradigm in the field of nonparametric Statistics and has been used in many applications, as illustrated in Sen [23, 24, 25, 26, 27, 28, 29, 30, 31] and Sen and Ghosh [32], among others. The related theory is available in many sources, among which we mention Serfling [34], Sen and Ghosh [32], Lee [14] or Sen and Singer [33].

The test is derived under the assumption that $\mathbb{E}[e_{ij}^4] < \infty$ and thus, accommodates a large class of distributions (not necessarily continuous) as the sources of variation in model (1.1).

In Section 2, using a martingale property, we describe the decomposition of the $U$-statistic that provides the key to the proposed test. In Section 3, we present simulation studies designed to evaluate the efficiency of the proposed test in different situations and to compare its performance with that of the exact $F$-test. We conclude with a brief discussion in Section 4.

## 2. Testing for variance components

Consider the function $g(x, y) = (x - y)^2/2$ and note that, under model (1.1),

$$\mathbb{E}[g(Y_{ij}, Y_{ij'})] = \mathbb{E}[(e_{ij} - e_{ij'})^2]/2 = \sigma_e^2.$$

Therefore, an unbiased estimator of $\sigma_e^2$, based only on the $n_i$ observations obtained under the $i$-th treatment is given by the following $U$-statistic

$$U_i = \binom{n_i}{2}^{-1} \sum_{1 \leq j < j' \leq n_i} g(Y_{ij}, Y_{ij'}) = \binom{n_i}{2}^{-1} \sum_{1 \leq j < j' \leq n_i} (Y_{ij} - Y_{ij'})^2/2$$

$$(2.1) \qquad = (n_i - 1)^{-1} \sum_{j=1}^{n_i} (Y_{ij} - \overline{Y}_{i.})^2 = S_i^2, \quad i = 1, \ldots, k.$$

Since $\mathbb{E}[(b_i - b_{i'})(e_{ij} - e_{i'j'})] = 0$, it follows that $\mathbb{E}[g(Y_{ij}, Y_{i'j'})] = \{2\sigma_b^2 + 2\sigma_e^2\}/2 = \sigma_b^2 + \sigma_e^2$. Therefore, an unbiased estimator of $\sigma_b^2 + \sigma_e^2$, based only on the observations obtained under treatments $i$ and $i'$, is given by the following generalized $U$-statistic of order (1,1)

$$(2.2) \qquad U_{ii'} = (n_i n_{i'})^{-1} \sum_{j=1}^{n_i} \sum_{j'=1}^{n_{i'}} \frac{(Y_{ij} - Y_{i'j'})^2}{2}, \quad 1 \leq i < i' \leq k.$$

Letting $n = \sum_{i=1}^{k} n_i$, the lexicographically ordered observations,

$$Y_{11}, \ldots, Y_{1n_1}, Y_{21}, \ldots, Y_{2n_2}, \ldots, Y_{k1}, \ldots, Y_{kn_k},$$

may be re-expressed as

$$Y_1, \ldots, Y_{n_1}, Y_{n_1+1}, \ldots, Y_{n_1+n_2}, \ldots, Y_{n_1+\cdots+n_{k-1}+1}, \ldots, Y_n,$$

where the first $n_1$ observations relate to treatment 1, the next $n_2$, to treatment 2 and so on. The uniformly minimum variance unbiased estimator (UMVUE) of the variance of the observations is given by the $U$-statistic

$$U_n^0 = \binom{n}{2}^{-1} \sum_{1 \leq r < s \leq n} \frac{1}{2}(Y_r - Y_s)^2,$$

which can be rewritten as

$$
\begin{aligned}
U_n^0 &= \binom{n}{2}^{-1} \left\{ \sum_{i=1}^{k} \binom{n_i}{2} U_i + \sum_{1 \leq i < i' \leq k} n_i n_{i'} U_{ii'} \right\} \\
(2.3) \qquad &= \sum_{i=1}^{k} \frac{n_i(n_i-1)}{n(n-1)} U_i + 2 \sum_{1 \leq i < i' \leq k} \frac{n_i n_{i'}}{n(n-1)} U_{ii'},
\end{aligned}
$$

highlighting its nature of a linear combination of generalized $U$-statistics. The first and second terms in (2.3) correspond, respectively, to the within and between-treatments components. Since

$$n_i(n_i - 1)/\{n(n-1)\} = \{n_i/n\} - \{n_i(n - n_i)\}/\{n(n-1)\},$$

and $\sum_{i' \neq i} n_{i'} = n - n_i$, we can rewrite the first term in (2.3) as

$$
\begin{aligned}
\sum_{i=1}^{k} \frac{n_i(n_i-1)}{n(n-1)} U_i &= \sum_{i=1}^{k} \frac{n_i}{n} U_i - \sum_{i=1}^{k} \frac{n_i(n - n_i)}{n(n-1)} U_i \\
&= \sum_{i=1}^{k} \frac{n_i}{n} U_i - \sum_{1 \leq i < i' \leq k} \frac{n_i n_{i'}}{n(n-1)} U_i - \sum_{1 \leq i < i' \leq k} \frac{n_i n_{i'}}{n(n-1)} U_{i'},
\end{aligned}
$$

leading to the decomposition

$$(2.4) \qquad U_n^0 = W_n + B_n$$

where

$$W_n = \sum_{i=1}^{k} \frac{n_i}{n} U_i \quad \text{and} \quad B_n = \sum_{1 \leq i < i' \leq k} \frac{n_i n_{i'}}{n(n-1)} \{2 U_{ii'} - U_i - U_{i'}\}.$$

This includes the classical ANOVA decomposition as a special case. Now observe that

$$\begin{aligned}\mathbb{E}[B_n] &= \sum_{1\leq i<i'\leq k}\frac{n_i n_{i'}}{n(n-1)}\{2(\sigma_b^2+\sigma_e^2)-2\sigma_e^2\} = 2\sigma_b^2\sum_{1\leq i<i'\leq k}\frac{n_i n_{i'}}{n(n-1)}\\ &= \sigma_b^2\sum_{i=1}^{k}\sum_{i'\neq i}\frac{n_i n_{i'}}{n(n-1)} = \sigma_b^2\sum_{i=1}^{k}\frac{n_i(n-n_i)}{n(n-1)} = \sigma_b^2\left(\frac{n^2-\sum_{i=1}^{k}n_i^2}{n(n-1)}\right)\geq 0.\end{aligned}$$

Therefore, $\mathbb{E}[B_n]=0$ if and only if $\sigma_b^2=0$. These results allow us to construct a test for (1.2) based on $B_n$ (suitably standardized) as we show in the sequel.

Firstly using the lexicographically ordered observations, we show that $B_n$ may be reexpressed as

$$(2.5) \qquad B_n = \binom{n}{2}^{-1}\sum_{1\leq r<s\leq n}\eta_{nrs}\psi(Y_r,Y_s),$$

where $\psi(x_1,x_2)=(x_1-\mu)(x_2-\mu)$ and

$$(2.6) \quad \eta_{nrs} = \begin{cases}(n-n_i)/(n_i-1), & \text{if } Y_r \text{ and } Y_s \text{ are both observed under the } i\text{-th treatment,}\\ -1, & \text{otherwise.}\end{cases}$$

Thus, it follows that

$$(2.7) \qquad \sum_{1\leq r<s\leq n}\eta_{nrs}=0$$

and

$$(2.8) \qquad \sum_{1\leq r<s\leq n}\eta_{nrs}^2 = \binom{n}{2}(k-1)\left\{1+\frac{1}{n}\sum_{i=1}^{k}\frac{n-n_i}{(n_i-1)(k-1)}\right\}.$$

When $k\to\infty$ and the $n_i$'s are bounded, i.e., the maximum number of observations per treatment does not increase with the number of treatments, it follows that $k=O(n)$; hence, by (2.8) we may conclude that $M_n=\sum_{1\leq r<s\leq n}\eta_{nrs}^2=O(n^3)$.

Under $\mathcal{H}_0:\sigma_b^2=0$, the random variables $Y_1,\ldots,Y_n$ are independent and identically distributed and $\psi(x,y)$ is a first-order stationary kernel, centered at 0, constituting an *orthogonal system* in the sense adopted in Pinheiro et al. [19], that is,

$$(2.9) \qquad \psi_1(Y_r) = \mathbb{E}[\psi(Y_r,Y_s)|Y_r] = \mathbb{E}[(Y_r-\mu)(Y_s-\mu)|Y_r] = 0 \text{ a.s.}$$
$$(2.10) \qquad \mathbb{E}[\psi(Y_1,Y_2)\psi(Y_1,Y_3)] = 0.$$

Now, observe that $\mathbb{E}[\psi^2(Y_r,Y_s)]=\text{Var}[\psi(Y_r,Y_s)]=\mathbb{E}[(Y_r-\mu)^2(Y_s-\mu)^2]=\sigma_e^4<\infty$. Then, we may show (see Appendix), that under $\mathcal{H}_0:\sigma_b^2=0$,

$$(2.11) \qquad J_n = \frac{\binom{n}{2}B_n}{W_n\sqrt{M_n}} \xrightarrow{D} \mathcal{N}(0,1)$$

when $k\to\infty(\Rightarrow n\to\infty)$. The key to the proof is the martingale property exhibited by $B_n$, as pointed in Pinheiro et al. [19] in a different setup. Let $\lim_{n\to\infty}M_n/n^3 =$ $\lambda$, which is finite because $M_n = O(n^3)$; assume that the fourth moment of the distribution of the random effects is finite and let $\delta$ denote a constant. Then, under the sequence of local hypotheses

$$\mathcal{H}_{1n} : \sigma_b^2 = \delta^2/\sqrt{n}, \tag{2.12}$$

it follows that

$$J_n \xrightarrow{D} \mathcal{N}\left(\frac{\delta^2}{2\sigma_e^2\sqrt{\lambda}}, 1\right) \text{ as } k \to \infty, \tag{2.13}$$

Details are presented in the Appendix. Because the mean of this limiting normal distribution is positive, we may use $J_n$ as a test statistic for (1.2), rejecting the null hypothesis $\mathcal{H}_0$ at significance level $\alpha$ when $J_n \geq z_\alpha$, where $z_\alpha$ is the $(1-\alpha)100\%$ percentile of the standard normal distribution. By (2.13), the power of this test is directly related to the magnitude of the intraclass correlation coefficient $\rho = \sigma_b^2/(\sigma_b^2 + \sigma_e^2)$; more specifically, the power is a monotone increasing function of $\rho$, as expected.

## 3. Simulation results

We summarize some simulation studies conducted with the objective of evaluating the behaviour of the proposed test. First, we examine the efficiency of the $U$-test in balanced studies under different distributions for the within-treatment errors $e_{ij}$ and random effects $b_i$. Then, we compare the efficiency of the proposed $U$-test with that of the usual $F$-test under various settings. Additional results from the simulation studies may be obtained from the authors.

### 3.1. Efficiency of the proposed test

To evaluate the behaviour of the proposed test for small and moderate samples we considered $10,000$ Monte-Carlo samples obtained under model (1.1) with $\mu = 2$ and $\sigma_e^2 = 1$, for different numbers of treatments ($k = 10, 30$ and $100$) in balanced studies. We assumed that $e_{ij} \sim \mathcal{N}(0,1)$ and $b_i \sim t_3 \times \sqrt{\sigma_b^2/3}$, so that $\text{Var}[b_i] = \sigma_b^2$. The between-treatments variance, $\sigma_b^2$, was set to 0 (to estimate the size of the test), 0.2, 0.5 and 1. The empirical power of the test under each setting was evaluated for a significance level of $\alpha = 0.05$.

The results, displayed in first half of Table 1, suggest that the $U$-test tends to be liberal when there are few observations for each treatment ($n_i \leq 4$). The difference between the empirical and nominal size of the test is acceptable when the number of treatments is at least 30 and there are 5 or more observations per treatment. Furthermore, the power of the test increases with $\sigma_b^2$, $n_i$ and $k$, as expected.

To illustrate the robustness of the $U$-test with respect to heavy-tailed within-treatment error distributions, we repeated the simulation study described previously assuming that $e_{ij} \sim t_5 \times \sqrt{3/5}$ (so that $\text{Var}[e_{ij}] = 1$). The results, summarized in the second half of Table 1, are not considerably different from the previous ones; here, however, the empirical size of the test is closer to the nominal size and the test also tends to be more powerful than when the within-treatment errors are normal.

Table 1

Rejection rates (%) for the 5% level $U$-test in balanced designs with $b_i \sim t_3 \times \sqrt{\sigma_b^2/3}$

| | $k = 10$ | | | | $k = 30$ | | | | $k = 100$ | | | |
|---|---|---|---|---|---|---|---|---|---|---|---|---|
| | $n_i$ | | | | $n_i$ | | | | $n_i$ | | | |
| $\sigma_b^2$ | 2 | 4 | 5 | 10 | 2 | 4 | 5 | 10 | 2 | 4 | 5 | 10 |
| | | | | | $e_{ij} \sim$ normal | | | | | | | |
| 0   | 14 | 8.8 | 7.6 | 6.2 | 8.3 | 7.1 | 6.3 | 6.1 | 7.2 | 6.1 | 5.9 | 5.4 |
| 0.2 | 26 | 36  | 41  | 64  | 33  | 59  | 69  | 93  | 54  | 92  | 97  | 100 |
| 0.5 | 41 | 61  | 69  | 89  | 61  | 91  | 96  | 100 | 93  | 100 | 100 | 100 |
| 1.0 | 57 | 82  | 88  | 98  | 85  | 99  | 100 | 100 | 100 | 100 | 100 | 100 |
| | | | | | $e_{ij} \sim t_5 \times \sqrt{3/5}$ | | | | | | | |
| 0   | 13 | 8.9 | 7.8 | 6.3 | 8.4 | 6.8 | 6.3 | 6.1 | 7.3 | 5.9 | 5.5 | 5.4 |
| 0.2 | 26 | 37  | 42  | 64  | 34  | 61  | 69  | 93  | 54  | 92  | 97  | 100 |
| 0.5 | 41 | 62  | 71  | 89  | 63  | 91  | 95  | 100 | 92  | 100 | 100 | 100 |
| 1.0 | 60 | 83  | 89  | 97  | 86  | 99  | 100 | 100 | 100 | 100 | 100 | 100 |

### 3.2. Comparison between the F- and U-tests

To compare the two tests, we first repeated the simulation study described previously, using the same values for $\mu$ and $\sigma_e^2$ in the following situations:

(i) Balanced study under normality;
(ii) Moderately unbalanced study under normality;
(iii) Lightly unbalanced study with heavy-tailed distributions;
(iv) Balanced study with heavy-tailed asymmetric distributions.

Firstly, we consider a balanced setting with 5 observations per treatment. The random effects and the within-treatment errors were generated according to: $b_i \sim \mathcal{N}(0, \sigma_b^2)$ and $e_{ij} \sim \mathcal{N}(0, 1)$. The results are presented in Table 2 and as in the previous studies, they suggest that the $U$-test is liberal. However, the difference between the empirical and nominal size of the test decreases substantially when the number of treatments increases. This occurs mainly when there are few units per treatment. On the other hand, the proposed $U$-test is more powerful than the exact $F$-test, though we should recognize that this might be influenced by its liberal nature.

We also considered simulation studies to evaluate the effect of imbalance in the power of the two tests under investigation. Khuri et al. [13], Donner and Koval [7] and Lee [15] discuss the performance of $F$-tests for different degrees of imbalance, defined as $\kappa = 1/(1 + \psi^2)$ where $\psi$ denotes the coefficient of variation associated with the sample sizes $n_1, \ldots, n_k$. By definition, $0 < \kappa \leq 1$ with $\kappa = 1$ only for balanced studies. Smaller values of $\kappa$ indicate larger degrees of imbalance.

We selected the number of within-treatment observations from a geometric distribution (shifted to the right by 2) with parameter $p = 0.15$ and support $\{0, 1, \ldots\}$. This corresponds to a moderate degree of imbalance ($\kappa = 0.61$). The random effects $b_i$ and the within-treatment errors $e_{ij}$ were respectively generated from $\mathcal{N}(0, \sigma_b^2)$ and $\mathcal{N}(0, 1)$ distributions.

In spite of the moderate imbalance, the size of the $F$-test is closer to the nominal level than is that of the (more liberal) $U$-test, especially when the number of treatments is small; however, the latter is more powerful than the $F$-test in most settings.

We also evaluated the simultaneous effects of a small degree of imbalance and heavier-tailed distributions for the random effects and within-treatment errors. In this direction, we selected the number of within-treatment observations with equal probability from the set $\{5, 6, 7, 8, 9, 10\}$; this corresponds to $\kappa \cong 0.95$. In addition, we assumed that $b_i \sim t_{4.1} \times \sqrt{\sigma_b^2 21/41}$ and $e_{ij} \sim t_{4.1} \times \sqrt{21/41}$, so that $\text{Var}[b_i] = \sigma_b^2$ and $\text{Var}[e_{ij}] = 1$.

TABLE 2
*Rejection rates (%) for the 5% level F and U-tests*

| $\sigma_b^2$ | $k=10$ | | $k=20$ | | $k=30$ | | $k=50$ | | $k=100$ | |
|---|---|---|---|---|---|---|---|---|---|---|
| | $F$ | $J_n$ | $F$ | $J_n$ | $F$ | $J_n$ | $F$ | $J_n$ | $F$ | $J_n$ |
| | Balanced designs and normal distributions | | | | | | | | | |
| 0 | 4.9 | 7.7 | 4.7 | 7.3 | 5.3 | 7.0 | 5.0 | 6.0 | 4.9 | 5.8 |
| 0.2 | 41 | 51 | 63 | 71 | 78 | 82 | 92 | 94 | 99 | 100 |
| 0.5 | 78 | 84 | 96 | 97 | 99 | 100 | 100 | 100 | 100 | 100 |
| 1 | 95 | 97 | 100 | 100 | 100 | 100 | 100 | 100 | 100 | 100 |
| | Moderately unbalanced designs and normal distributions | | | | | | | | | |
| 0 | 5.2 | 7.2 | 5.0 | 6.7 | 5.1 | 6.2 | 4.9 | 6.0 | 4.9 | 5.7 |
| 0.2 | 51 | 70 | 74 | 87 | 90 | 96 | 98 | 100 | 100 | 100 |
| 0.5 | 96 | 99 | 100 | 100 | 100 | 100 | 100 | 100 | 100 | 100 |
| 1 | 100 | 100 | 100 | 100 | 100 | 100 | 100 | 100 | 100 | 100 |
| | Unbalanced designs with heavy-tailed distributions for $b_i$ and $e_{ij}$ | | | | | | | | | |
| 0 | 3.5 | 5.9 | 3.9 | 5.5 | 4.4 | 5.4 | 4.7 | 5.3 | 5.3 | 5.1 |
| 0.2 | 47 | 55 | 61 | 72 | 71 | 80 | 82 | 89 | 92 | 95 |
| 0.5 | 68 | 79 | 85 | 90 | 91 | 94 | 96 | 98 | 98 | 99 |
| 1 | 84 | 90 | 94 | 97 | 97 | 98 | 98 | 99 | 99 | 99 |
| | Balanced designs with heavy-tailed asymmetric distributions for $b_i$ and $e_{ij}$ | | | | | | | | | |
| 0 | 4.4 | 5.7 | 4.4 | 5.6 | 4.3 | 5.4 | 4.4 | 5.3 | 4.6 | 5.2 |
| 0.2 | 37 | 45 | 57 | 64 | 70 | 76 | 86 | 89 | 98 | 98 |
| 0.5 | 70 | 77 | 89 | 91 | 96 | 97 | 99 | 99 | 100 | 100 |
| 1 | 87 | 91 | 98 | 98 | 99 | 100 | 100 | 100 | 100 | 100 |

The results, also presented in Table 2, suggest that the size of the $U$-test is closer to the nominal level than is that of the $F$-test even in settings with few treatments. The bias is larger for the $F$-test, mainly when there are few treatments.

Finally, we evaluated the possible effect of asymmetric heavy-tailed distributions for the random effects and within-treatment errors on the size and power of the two tests. In this direction, we considered a balanced setting with 5 observations per treatment. The random effects and the within-treatment errors were generated according to the following distributions: $b_i \sim \{(Y_1 - \mathbb{E}[Y_1])/\sqrt{\text{Var}[Y_1]}\} \times \sigma_b$ and $e_{ij} \sim (Y_2 - \mathbb{E}[Y_2])/\sqrt{\text{Var}[Y_2]}$, where $Y_1, Y_2$ are independent identically distributed random variables with skew $t$ distribution with 4.1 degrees of freedom, location parameter 0, dispersion parameter 1 and asymmetry parameter $\lambda = 1$ with index of skewness equal to 1.77. For details on the skew $t$ distribution, see Azzalini and Capitanio [2]. The results are also presented in Table 2 and suggest a conservative behaviour for the $F$-test. On the other hand, the size of the $U$-test is closer to the nominal level, even for samples of moderate size. It may also be observed that the $U$-test is more powerful in all configurations.

## 4. Discussion and conclusion

Although there exists an exact $F$-test with optimal properties for testing the significance of the between-treatments variance component in a one-way random effects model with balanced data under normality, we must rely on sub-optimal or approx-imate tests in unbalanced or nonnormal settings. We derived an asymptotic $U$-test that may be employed with unbalanced data and does not require a specified form for the underlying probability distributions. Simulation studies suggest that it has reasonable properties even for moderate and small samples.

In particular, we may conclude that the $F$-test is more affected by the lack of normality of the random effects and within-treatment errors than by imbalance. Furthermore, the $U$-test seems to be less sensitive to imbalance and to be more powerful than the $F$-test in general. Such conclusions must be viewed with caution, given the liberal nature of the $U$-test, specially for small sample sizes. To bypass this problem, one could consider bootstrap methods to obtain the empirical distribution of the statistic $J_n$ under $\mathcal{H}_0$ and use the fact that the null hypothesis should be rejected for high values of $J_n$. For details on the use of bootstrap to construct confidence regions or for hypothesis testing, see Davison and Hinkley [5] or Lehmann and Romano [16], Chapter 15, for example.

In summary, under normality or under moderate degrees of imbalance, the $F$-test behaves well in terms of power when compared to $U$-test; however, in situations where the distribution of the random effects and within-treatment errors are nonnormal, the $U$-test is preferable even when the number of treatments is small. In all settings the $U$-test is more powerful than the $F$-test, mainly for small and moderate samples.

The proposed $U$-test may also be used to test for simultaneous significance of all random components in models with more than one random effect, since it behaves as (2.11) under $\mathcal{H}_0$. The corresponding power computations, however, require further research.

Brownie and Boos [4] and Boos and Brownie [3] study the behaviour of $F$-tests and analogous rank tests in one-way random effects models when the number of treatments is large. Akritas and Arnold [1] also study the asymptotics of ANOVA under similar conditions but their results require more restrictive assumptions than those we considered; the proofs, however, are different for balanced and unbalanced studies.

In the context under investigation, the derivation of tests for (1.2) may not follow the standard procedures since the null hypothesis defines a point (or region) on the boundary of the parameter space and this brings in some technical difficulties. Asymptotic tests for (1.2) or, more generally, for testing the significance of variance components under linear mixed models are available in the literature. Based on the ideas of Silvapulle and Silvapulle [35], Verbeke and Molenberghs [38] obtained score-type tests under the assumption that the underlying probability distributions are normal. Along the same lines, Savalli et al. [20] extended the results to accommodate elliptical underlying distributions. In particular, for the one-way random effects models, the corresponding test statistic follows an asymptotic distribution given by a 50:50 mixture of $\chi_0^2$ and $\chi_1^2$ distributions. Tests based on generalized likelihood methods (that are asymptotically equivalent to the score-type tests) are considered in Self and Liang [22], Stram and Lee [37] and Silvapulle and Sen [36], for example. The main disadvantage of such tests is the difficulty in verifying the required regularity conditions as shown in Giampaoli and Singer [9]. The derivation of the proposed $U$-test is not affected by such difficulties and we envisage that it may serve as a building block for more general setup as indicated in Nobre [18].

Other alternatives have been suggested in the literature. Using a Laplace expansion of the integrated log-likelihood, Lin [17] obtained a global score test for the hypothesis that all variance components are zero in the framework of generalized linear models with random effects. Along the same lines, Zhu and Fung [39] ob-tained a global score-type test for the hypothesis that all variance components are zero in a semiparametric mixed model, under the assumption of normality only of the conditional error vector. It is important to point out, that the tests are obtained without specifying the distribution of the random effects, though they are based on a marginal approach under which the within-treatment covariances may be negative. Hall and Præstgaard [10] consider such a constraint in the framework of generalized linear models with random effects, via a projected score test. In practice, their results are difficult to apply when the dimension of the vector of random effects is large.

**Appendix A: Proof of (2.11)**

Note that under $\mathcal{H}_0 : \sigma_b^2 = 0$,

$$B_n = \binom{n}{2}^{-1} \sum_{1 \leq r < s \leq n} \eta_{nrs} \psi(Y_r, Y_s)$$

is a sum of $\binom{n}{2}$ uncorrelated terms such that $\mathbb{E}[\psi(Y_r, Y_s)] = 0$ e $\text{Var}[\psi(Y_r, Y_s)] = \sigma_e^4 < \infty$. Thus, it follows that

(A.1) $$\mathbb{E}[B_n] = 0 \quad \text{and} \quad \text{Var}[B_n] = \sigma_e^4 M_n \binom{n}{2}^{-2}.$$

Then, along the lines adopted in Pinheiro et al. [19], we explore the martingale structure of $B_n$, and apply a martingale Central Limit Theorem as given in Dvoretzky [8] or Sen and Singer [33] to obtain the desired asymptotic distribution. Their proof requires $k$ fixed and $n_i \to \infty$; we consider $n_i$ bounded and let $k \to \infty$, although the test is also valid when both $k \to \infty$ and $n_i \to \infty$.

Initially, consider the statistic

(A.2) $$T_n = \sum_{1 \leq i < j \leq n} \eta_{nij} \phi(X_i, X_j)$$

with $\phi$ satisfying

(A.3) $\phi_1(X_1) = \mathbb{E}[\phi(X_1, X_2) | X_1] = 0$ a.s.,
(A.4) $\mathbb{E}[\phi(X_1, X_2) \phi(X_1, X_3)] = 0,$
(A.5) $\mathbb{E}[\phi^2(X_1, X_2)] < \infty,$

where $X_1, \ldots, X_n$ represents a sequence of independent and identically distributed random variables and $\eta_{nij}$ denotes weights such that

(A.6) $$\sum_{1 \leq i < j \leq n} \eta_{nij} = 0 \quad \text{and} \quad \sum_{1 \leq i < j \leq n} \eta_{nij}^2 = M_n^*.$$

Defining, $Z_{nj} = \sum_{i=1}^{j-1} \eta_{nij} \phi(X_i, X_j)$, for $j = 2, \ldots n$ and $T_{nu} = \sum_{l=2}^{u} Z_{nl}$ for $2 \leq u \leq n$, it follows that $T_{nn} = T_n$. Now, considering the nondecreasing sequence (in $u$) of $\sigma$-fields $\sigma_{nu} = \sigma(X_i, i \leq u)$ that corresponds to the $\sigma$-fields generated by the random vector $(X_1, \ldots, X_u)^\top$, $2 \leq u \leq n$ we may apply the following result due to Pinheiro et al. [19].

**Theorem A.1.** *Under conditions (A.2)–(A.5), $\{T_{nu}, 2 \leq u \leq n\}$ is a zero mean martingale process adapted to the filter $\{\sigma_{nu}, 2 \leq u \leq n\}$.*

Letting $D_n$ denote a weakly consistent estimator of $\tau^2 = \mathbb{E}[\phi^2(X_i, X_j)] < \infty$ ($i < j$), it follows from (A.4) and (A.5) that $\text{Var}[T_n] = \mathbb{E}[T_n^2] = \tau^2 M_n^*$. Using Slutsky's theorem, we obtain

$$\text{(A.7)} \qquad \frac{D_n M_n^*}{\mathbb{E}[T_n^2]} \xrightarrow{\mathbb{P}} 1.$$

Now, note that the statistic $\binom{n}{2} B_n$ is of the form (A.2) with $\phi(x,y) = (x-\mu)(y-\mu)$ satisfying (A.3)–(A.5). Additionally, observe that $\tau^2 = \mathbb{E}[\phi^2(X_i, X_j)] = \sigma_e^4 < \infty$. A consistent estimator of $\tau^2$ may be obtained from the following result.

**Lemma A.1.** *Consider a linear combination of $U$-statistics of the form $W_n = \sum_{i=1}^{k} \frac{n_i}{n} U_i = \sum_{i=1}^{k} w_i U_i$. If*

$$\text{(A.8)} \qquad \max_{1 \leq i \leq k} \frac{n_i}{\sum_{i=1}^{k} n_i} \to 0,$$

*it follows that $W_n^2 \xrightarrow{\mathbb{P}} \sigma_e^4$.*

*Proof.* Condition (A.8) is valid when the $n_i$'s are fixed and $k \to \infty$. It also holds when $n_i \to \infty$, as for example, in a balanced study.

Noting that $\mathbb{E}[W_n] = \sigma_e^2$, and using properties of $U$-statistics, it is possible to show that $\text{Var}[U_i] = \mathbb{E}[e_{ij}^4]/n_i - (n_i - 3)\sigma_e^4/\{(n_i - 1)n_i\}$. Then, a direct application of Chebyshev's inequality leads to $W_n \xrightarrow{\mathbb{P}} \sigma_e^2$. Since $x^2$ is a continuous function in the support of $W_n$, we conclude that $D_n = W_n^2 \xrightarrow{\mathbb{P}} \sigma_e^4$.

In the proof, $\mathbb{E}[e_{ij}^4] < \infty$ is used to guarantee the existence of $\text{Var}[U_i]$. In some situations, as in the balanced case, this requirement may be relaxed and the result follows from Khintchine's weak law of large numbers. □

We now focus on the following lemma.

**Lemma A.2.** *Consider the weights defined in (2.6). When $k \to \infty$ it follows that*

$$\text{(A.9)} \qquad \sum_{1 \leq i \neq j < u \leq n} \eta_{niu}^2 \eta_{nju}^2 / M_n^2 \to 0, \quad \text{and} \quad \sum_{1 \leq i \neq j \leq n} \eta_{nij}^4 / M_n^2 \to 0.$$

*Proof.* When $k \to \infty$ and the $n_i$'s are bounded, it follows that $k = O(n)$. Then by (2.8), we obtain $M_n = \sum_{1 \leq r < s \leq n} \eta_{nrs}^2 = O(n^3)$. Now, the result follows from

$$\sum_{1 \leq i \neq j \leq n} \eta_{nij}^4 = O(n^5) \quad \text{and} \quad \sum_{1 \leq i \neq j < u \leq n} \eta_{niu}^2 \eta_{nju}^2 = O(n^5).$$

□

The proof of the main result consists in verifying the two regularity conditions of Theorem 3.3.7 in (Sen and Singer [33], page 120) and it is similar to the proof presented in Pinheiro et al. [19].

**Theorem A.2.** *Suppose that $\mathbb{E}[\phi^4(X_1, X_2)] < \infty$, (A.3)-(A.5) hold and, for $n \to \infty$, the weights $\eta_{nij}$ satisfy*

$$\text{(A.10)} \qquad \sum_{1 \leq i \neq j < k \leq n} \eta_{nik}^2 \eta_{njk}^2 / (M_n^*)^2 \to 0 \quad \text{and} \quad \sum_{1 \leq i \neq j \leq n} \eta_{nij}^4 / (M_n^*)^2 \to 0.$$

Then, for $T_n$ in (A.2), we have that

(A.11) $$(M_n^* D_n)^{-1/2} T_n \xrightarrow{D} \mathcal{N}(0,1).$$

*Proof.* For the weights defined in (2.6), we have $M_n^* = M_n$. Then (2.11) follows from Lemmas A.1 and A.2. □

## Appendix B: Proof of (2.13)

First observe that

$$U_i = \binom{n_i}{2}^{-1} \sum_{1 \le j < j' \le n_i} \frac{(Y_{ij} - Y_{ij'})^2}{2} = \binom{n_i}{2}^{-1} \sum_{1 \le j < j' \le n_i} \frac{(e_{ij} - e_{ij'})^2}{2}$$

and $W_n = \sum_{i=1}^{k} (n_i/n) U_i$, are functions only of the vector of within-treatment errors $\mathbf{e}$. Consequently, their distributions under $\mathcal{H}_0$ and $\mathcal{H}_{1n}$ are the same. On the other hand, we have

$$\begin{aligned}
B_n &= \sum_{1 \le i < i' \le k} \frac{n_i n_{i'}}{n(n-1)} \{2U_{ii'} - U_i - U_{i'}\} \\
&= \sum_{1 \le i < i' \le k} \frac{n_i n_{i'}}{n(n-1)} \left\{ \frac{2}{n_i n_{i'}} \sum_{j=1}^{n_i} \sum_{j'=1}^{n_{i'}} \frac{(Y_{ij} - Y_{i'j'})^2}{2} - U_i - U_{i'} \right\} \\
&= B_n^0 + C_n,
\end{aligned}$$

where

$$B_n^0 = \sum_{1 \le i < i' \le k} \frac{n_i n_{i'}}{n(n-1)} \left\{ \frac{2}{n_i n_{i'}} \sum_{j=1}^{n_i} \sum_{j=1}^{n_{i'}} \frac{(e_{ij} - e_{i'j'})^2}{2} - U_i - U_{i'} \right\},$$

$$C_n = \sum_{1 \le i < i' \le k} \frac{n_i n_{i'}}{n(n-1)} \left\{ (b_i - b_{i'})^2 + \frac{2(b_i - b_{i'})}{n_i n_{i'}} \sum_{j=1}^{n_i} \sum_{j'=1}^{n_{i'}} (e_{ij} - e_{i'j'}) \right\}.$$

Also, under $\mathcal{H}_0$, $B_n = B_n^0$. Therefore,

(B.1) $$J_n^* = \frac{\binom{n}{2} B_n}{\sigma_e^2 M_n^{1/2}} = \frac{\binom{n}{2} B_n^0}{\sigma_e^2 M_n^{1/2}} + \frac{\binom{n}{2} C_n}{\sigma_e^2 M_n^{1/2}} = J_{n0} + Q_n,$$

where $J_{n0} = \binom{n}{2} B_n^0 / (\sigma_e^2 M_n^{1/2}) \xrightarrow{D} \mathcal{N}(0,1)$ (by 2.11) and $Q_n = \binom{n}{2} C_n / (\sigma_e^2 M_n^{1/2})$. The term $C_n$ can be decomposed as $C_n = C_{n1} + C_{n2}$ where

$$C_{n1} = \frac{1}{n(n-1)} \sum_{1 \le i \le i' \le k} n_i n_{i'} (b_i - b_{i'})^2$$

and

$$C_{n2} = \frac{2}{n(n-1)} \sum_{1 \le i \le i' \le k} (b_i - b_{i'}) \sum_{j=1}^{n_i} \sum_{j'=1}^{n_{i'}} (e_{ij} - e_{ij'}).$$

Since $\mathbb{E}[C_{n2}] = 0$, we have

$$\mathbb{E}[C_n] = \mathbb{E}[C_{n1}] = \frac{1}{n(n-1)} \sum_{1 \leq i \leq i' \leq k} n_i n_{i'} \mathbb{E}[(b_i - b_{i'})^2]$$

$$= \frac{1}{n(n-1)} \sum_{1 \leq i \leq i' \leq k} n_i n_{i'} 2\sigma_b^2 = \sigma_b^2 \left( \frac{n^2 - \sum_{i=1}^k n_i^2}{n(n-1)} \right) = O(n^{-1/2}).$$

The following result provides the asymptotic behaviour of $Q_n$ under the assumption of the existence of the fourth moment of the random effects.

**Lemma B.1.** *Consider model (1.1) under $\mathcal{H}_{1n} = \sigma_b^2 = \delta^2/\sqrt{n}$. Assuming that $\mathbb{E}[b_i^4] < \infty$ and letting $\lim_{n \to \infty} M_n/n^3 = \lambda$, we have*

$$\frac{\binom{n}{2} C_n}{\sqrt{M_n}} \xrightarrow{\mathbb{P}} \frac{\delta^2}{2\sqrt{\lambda}}.$$

*Proof.* Given that $\mathbb{E}[b_i^4] < \infty$, then by contiguity of the sequence of local alternative hypotheses $\mathcal{H}_{1n} : \sigma_b^2 = \delta^2/\sqrt{n}$, we have that $\mathbb{E}[b_i^4] = o(1)$ under $\mathcal{H}_{1n}$. Defining $b^{ij} = (b_i - b_j)^2$, recalling that $b_1, \ldots, b_k$ are independent and identically distributed and that the $n_i$'s are bounded, we have

$$\text{Var}[C_{n1}] = \frac{1}{n^2(n-1)^2} \left\{ \sum_{1 \leq i < i' \leq k} n_i^2 n_{i'}^2 \text{Var}[b^{ii'}] + \sum_{\substack{1 \leq i < t \leq k \\ 1 \leq i < w \leq k \\ t \neq w}} n_i^2 n_t n_w \text{Cov}[b^{it}, b^{iw}] \right\}$$

$$\leq \frac{\max\{\text{Var}[b^{12}], \text{Cov}[b^{12}, b^{13}]\}}{n^2(n-1)^2} \left\{ \sum_{1 \leq i < i' \leq k} n_i^2 n_{i'}^2 + \sum_{\substack{1 \leq i < t \leq k \\ 1 \leq i < w \leq k \\ t \neq w}} n_i^2 n_t n_w \right\}$$

$$= o(1) O(n^{-1}) = o(n^{-1}),$$

since $\text{Var}[b_a^{12}]$ and $\text{Cov}[b_a^{12}, b_a^{13}]$ are of the same order as $\mathbb{E}[b_i^4]$, namely, $o(1)$. Similarly, it is possible to show that $\text{Var}[C_{n2}] = o(n^{-1})$. Here, observing that $\binom{n}{2}/\sqrt{M_n} = O(\sqrt{n})$, we have

$$\lim_{n \to \infty} \frac{\binom{n}{2} \mathbb{E}[C_{n1}]}{\sqrt{M_n}} = \lim_{n \to \infty} \frac{\delta^2 \left( n^2 - \sum_{i=1}^k n_i^2 \right)}{2\sqrt{n}\sqrt{M_n}} = \frac{\delta^2}{2\sqrt{\lambda}},$$

and since $O(n) \text{Var}[C_{n1}] = o(1) \to 0$, when $n \to \infty$, it follows that $\binom{n}{2} C_{n1}/\sqrt{M_n} \xrightarrow{\mathbb{P}} \delta^2/\{2\sqrt{\lambda}\}$. Since $\mathbb{E}[C_{n2}] = 0$ and $\text{Var}[C_{n2}] = o(n^{-1})$, we have that $\binom{n}{2} C_{n2}/\sqrt{M_n} \xrightarrow{\mathbb{P}} 0$. Therefore,

$$\binom{n}{2} C_n/\sqrt{M_n} = \binom{n}{2}(C_{n1} + C_{n2})/\sqrt{M_n} \xrightarrow{\mathbb{P}} \delta^2/\{2\sqrt{\lambda}\}.$$

□

Now, recalling (B.1), Lemma B.1 observing that $W_n \xrightarrow{\mathbb{P}} \sigma_e^2$, and using Slutsky's theorem, we obtain

$$J_n = \frac{\binom{n}{2} B_n}{W_n \sqrt{M_n}} = \frac{\binom{n}{2} B_n^0}{\sigma_e^2 \sqrt{M_n}} \frac{\sigma_e^2}{W_n} + \frac{\binom{n}{2} C_n}{\sigma_e^2 \sqrt{M_n}} \frac{\sigma_e^2}{W_n} \xrightarrow{D} \mathcal{N}\left( \frac{\delta^2}{2\sigma_e^2 \sqrt{\lambda}}, 1 \right),$$

implying that, under $\mathcal{H}_{1n} : \sigma_b^2 = \delta^2/\sqrt{n}$ the center of the normal is shifted to the right by $\delta^2/(2\sigma_e^2\sqrt{\lambda})$.

**Acknowledgments.** We are grateful to Professor P. K. Sen and Drs. N. I. Tanaka and A. S. Pinheiro for constructive comments and suggestions. It is our pleasure to contribute to this Festschrift in the honour of Professor Sen.